\documentclass[11pt, twoside]{article}
\usepackage{amsfonts}
\usepackage{CJK,CJKnumb}
\usepackage{amssymb}
\usepackage{amsmath}
\usepackage{amsthm}
\usepackage{xcolor}
\usepackage{mathrsfs}
\usepackage[active]{srcltx}
\usepackage{mathrsfs}
\usepackage{makeidx}

\usepackage{amsfonts}
\usepackage{appendix}
\usepackage[active]{srcltx}
\usepackage{mathrsfs}
\usepackage{makeidx}

\allowdisplaybreaks

\allowdisplaybreaks

\renewcommand\hat{\widehat}

\def\supp{{\rm{\,supp\,}}}

\newtheorem{theorem}{Theorem}[section]
\newtheorem{lemma}[theorem]{Lemma}
\newtheorem{corollary}[theorem]{Corollary}
\newtheorem{proposition}[theorem]{Proposition}
\theoremstyle{definition}
\newtheorem{example}[theorem]{Example}
\newtheorem{remark}[theorem]{Remark}

\newtheorem{definition}[theorem]{Definition}
\numberwithin{equation}{section}

\numberwithin{equation}{section}

\numberwithin{equation}{section}

\begin{document}
\arraycolsep=1pt
\title{\Large\bf Notes on a Special Order on $\mathbb{Z}^\infty$}
\author{Jiawei Sun, Chao Zu\thanks{Corresponding author.} \ and Yufeng Lu}
\maketitle
\vspace{-0.8cm}
\section{Abstract}
\quad In 1958, Helson and Lowdenslager extended the theory of analytic functions to a general class of groups with ordered duals. In this context, analytic functions on such a group $G$ are defined as the integrable functions whose Fourier coefficients lie in the positive semigroup of the dual of $G$. In this paper, we found some applications of their theory to infinite-dimensional complex analysis. Specifically, we considered a special order on $\mathbb{Z}^\infty$ and corresponding analytic continuous functions on $\mathbb{T}^\omega$, which serves as the counterpart of the disk algebra in infinitely many variables setting. By characterizing its maximal ideals, we have generalized the following theorem to the infinite-dimensional case:
For a positive function $w$ that is integrable and log-integrable on $\mathbb{T}^d$, there exists an outer function $g$ such that $w=|g|^2$ if and only if the support of $\hat{\log w}$ is a subset of $\mathbb{N}^d\cap (-\mathbb{N})^d$. Furthermore, we have found the counterpart of the above function algebra in the closed right half-plane, and the representing measures of each point in the right half-plane for this algebra. As an application of the order, we provided a new proof of the infinite-dimensional Szeg\"{o}'s theorem.

\section{Introduction}
\quad The theory of analytic functions on the unit disk has been generalized in two cases: one is replacing the unit disk with plane domains, and the other replacing it with ordered groups, i.e., a group equipped with a translation invariant total order \cite{hl}. More specifically, the first generalization of analytic functions is the class of functions on $\mathbb{T}^d$ (or $\mathbb{D}^d$) whose Fourier transforms are supported on the semigroup generated by the copies of $\mathbb{N}$. The main purpose of this paper is to extend some classical theorems on this type of Hardy spaces to the infinite-dimensional case; however, the methods we use involve another type of analytic functions whose definition is related to a special order on $\mathbb{Z}^\infty$.

The way to equip $\mathbb{Z}^\infty$ (the discrete group generated by countable copies of $\mathbb{Z}$) with a translation invariant total order is not unique, such as lexicographic order, i.e., $\{x_i\}_{i=1}^\infty\geq \{y_i\}_{i=1}^\infty$ if and only if $x_\alpha>y_\alpha$ for the first index $\alpha$ that $x_\alpha\neq y_\alpha$. In this paper we equip $\mathbb{Z}^\infty$ with an order by embedding $\mathbb{Z}^\infty$ into the real line $\mathbb{R}$. To be specific, $$\tau: \{x_i\}_{i=1}^\infty\to \sum_{i=1}^\infty x_i\log p_i$$ is a homomorphism from $\mathbb{Z}^\infty$ to a subgroup of $\mathbb{R}$, where $p_i$ is the $i$-th prime, thus we can restrict the order of $\mathbb{R}$ to $\mathbb{Z}^\infty$. This idea is inspired by the research of Defant and Schoolmann on Dirichlet series \cite{ds3}. Their work further develpoed Bohr's view on Dirichlet series \cite{b}. Specifically, let $p_i$ be the $i$-th prime, then each positive integer $k$ has a unique factorization:
$k=\prod_{j=1}^\infty p_j^{\alpha_j}$
. If we replace each $k$ by the corresponding sequence $\{\alpha_j\}$ of exponents, the Dirichlet series
$\sum_{k=1}^\infty c_k/k^s$ takes the form
\begin{align}\label{e2.2}
\sum_{\alpha_i\geq 0}c(\alpha_1,\alpha_2,...)e^{-s\sum_{j=1}^\infty \alpha_j \log p_i}.
\end{align}
Writing $z_j=p_j^{-s}$, this series becomes a power series in infinitely many variables, namely, 
\begin{align}\label{e2.1}
\sum_{\alpha_i\geq 0}c(\alpha_1,\alpha_2,...)z_1^{\alpha_1}z_2^{\alpha_2}...
\end{align}
If $s$ is pure imaginary, this is a trigonometric series on $\mathbb{T}^\omega$, the Cartesian product of countable copies of torus $\mathbb{T}$. The map from \eqref{e2.1} to \eqref{e2.2} is called the Bohr transform. There is a known result that the Bohr transform is an isometry form $H^\infty_+(\mathbb{T}^\omega)$ to $\mathcal{H}^\infty$ \cite[Theorem 3.8\& Theorem 5.1]{ddm}, where $H^\infty_+(\mathbb{T}^\omega)$ is the space of essential bounded functions on $\mathbb{T}^\omega$ that take the form of \eqref{e2.1}, and $\mathcal{H}^\infty$ is the space of Dirichlet series that are bounded on the right half-plane.

In some cases, the power series \eqref{e2.1} corresponds to a holomorphic function in the domains $B_{c_0}$ or $B_{c_0}\cap l^2$ \cite[chapter 2\& chapter 13]{ddm}, where $B_{c_0}$ is the unit ball under the supremum norm in the space formed by Cauchy sequences that converge to $0$, and the term "holomorphic" refers to Fréchet differentiable \cite[chapter 2, section 2]{ddm} at each point of the corresponding domain. For $p\in [1,\infty]$, we define $H^p_+(\mathbb{T}^\omega)$ as the space of $L^p(\mathbb{T}^\omega)$-functions, whose Fourier series take the form of \eqref{e2.1}. $H^p_+(\mathbb{T}^\omega)$ is isomorphic to $H^p_+(l^2\cap B_{c_0})$ for $1<p<\infty$ \cite[Theorem 13.2]{ddm}, where $H^p_+(B_{c_0}\cap l^2)$ is defined as 
\begin{align*}
\left\{f \text{ is holomorphic on $B_{c_0}\cap l^2$}:
\sup_{0<r<1, N\in \mathbb{N}}\|f(rw_1,...,rw_N,0,0...)\|_{L^p(\mathbb{T}^N)}<\infty\right\}.
\end{align*}
And when $p=\infty$, $H^p_+(\mathbb{T}^\omega)$ is isomorphic to $H^p_+(B_{c_0})$, where the definition of $H^p_+(B_{c_0})$ is similar to that of $H^p_+(B_{c_0}\cap l^2)$.
The space $H^p_+(\mathbb{T}^\omega)$ has received extensive study in recent years, partly because it is the most natural generalization of Hardy spaces in the theory of several complex variables, and partly due to its connection with Dirichlet series.

In recent years, Defant and Schoolmann extended Bohr's theory to Dirichlet groups, such a topological group can serve as a compactification of real line \cite{ds,ds2,ds3}. Our research was mainly inspired by a profound result they obtained \cite[Theorem 2.13]{ds2}. From a one-sided perspective, they define the Cesàro mean on infinitely many variables setting: 
For $\mathbf{n}=\{n_1, n_2,...,0,0,....\}\in \mathbb{N}^\infty$, the Cesàro mean $\sigma_{x
}(\sum c_\mathbf{n})$ of formal series $\sum_{\mathbf{n}\in \mathbb{N}^\infty} c_\mathbf{n}$ is defined as
$\sum_{\tau(\mathbf{n})<x}c_\mathbf{n}(1-\frac{\tau(\mathbf{n})}{x})$, where $x$ is in the positive half of the real axis. This summation can naturally be viewed as being carried out according to the order mentioned above. The result shows that the Cesàro mean of the Fourier series of $f\in H^p_+(\mathbb{T}^\omega),\ p\in [1,\infty)$ on $\mathbb{T}^\omega$ will converge to $f$ with respect to $L^p$-norm as $x\to \infty$. This gives us a conviction that this order is significant.

Fortunately, there are many well-established results regarding Fourier analysis on groups with ordered dual \cite{hl, rud}, in which context, analytic functions on a topological group $G$ with an ordered dual $\Gamma$ are defined as the functions whose Fourier coefficients lie in the positive semigroup consisting of $\gamma\in \Gamma$ with $\gamma>0$. The Hardy space $H^p(G)$ corresponding to ordered group $\Gamma$ is defined as the space of analytic functions, as we mentioned above, in $L^p(G)$. Compared to $H^p_+(\mathbb{T}^\omega)$, $H_\tau^p(\mathbb{T}^\omega)$ better inherits the Fourier characteristics of Hardy spaces on the unit disk. However, the price to pay is that it may be difficult to find a class of holomorphic functions (on the Cartesian product of copies of unit disk) corresponding to such analytic functions.

In this paper, we shall be concerned with two problems. The first is a theorem in the theory of several complex variables \cite[Lemma 4]{n}: For a function $w$ that is integrable and log-integrable on $\mathbb{T}^d$, there exists an outer function $g\in H^p_+(\mathbb{T}^d)$ such that $w=|g|^2$ if and only if the support of $\hat{\log w}$ is a subset of $\mathbb{N}^d\cap (-\mathbb{N})^d$. We aim to extend it to the infinite-dimensional case. The second is the characterization of the maximal ideal spaces of the disk algebra $A_\tau(\mathbb{T}^\omega)$ on $\mathbb{T}^\omega$ corresponding to $\tau$, i.e., the analytic (corresponding to $\tau$) continuous functions on $\mathbb{T}^\omega$.

We have obtained the following result: First, we characterized the maximal ideal space of $A_\tau(\mathbb{T}^\omega)$, found the representing measures of its elements on $\mathbb{T}^\omega$, and used this to extend the aforementioned theorem to $\mathbb{T}^\omega$. Inspired by Bohr's vision on Dirichlet series, we further found the counterparts of $A_\tau(\mathbb{T}^\omega)$ and its representing measures on the right half-plane. Furthermore, as an application of order $\tau$ and related Hardy space, we provide a new proof of Szeg\"{o}'s theorem in the setting of infinitely many variables. The main techniques we use are almost entirely derived from Fourier analysis rather than from infinite-dimensional complex function theory.

\section{Preliminary}
\quad In this section we first recall some definitions, notations and theorems that are used throughout the paper. Most of these definitions originate from Fourier analysis on locally compact abelian groups. Following the notational conventions in this field, throughout this paper, $G$ always represents a locally compact abelian group, its Pontryagin dual is denoted by $\Gamma$, and $m_G$ denotes the normalized Haar measure on $G$. In the interest of simplicity, we write $L^p(m_G)$ as $L^p(G)$.





\begin{definition}\label{d3.1}
We call an abelian group $\Gamma$ an {\it ordered group} if there is a fixed semigroup $\Gamma'\subset\Gamma$ such that $\Gamma'\cap (-\Gamma')=0, \Gamma'\cup(-\Gamma')=\Gamma$.
Such a semigroup induces an translate invariant order ($x>y$ is equivalent to $x+z\geq y+z$, $\forall z\in \Gamma$) on $\Gamma$: For $x,y\in \Gamma$,  $x\geq y$ if $x-y\in \Gamma'$.
\end{definition}

\begin{example}\label{x1.1}
We denote by $\mathbb{Z}^\infty$ the discrete group generated by countable copies of $\mathbb{Z}$.
Every $\mathbf{n}\in \mathbb{Z}^\infty$ has the form
$$\mathbf{n}=(n_1,...n_k,0,0,0,...),\ n_i \in \mathbb{Z}.$$
For convenience, throughout the paper we replace $(n_1,...,n_k,0,0,...)$ by $(n_1,...,n_k)$, and we call $k$ the length of $(n_1,...,n_k,0,0,...)$ if $n_k$ is the last nonzero coordinate.

There is an archimedean total order $\tau$ on $\mathbb{Z}^\infty$ induced by the semigroup $\{\mathbf{n}\in \mathbb{Z}^\infty: \sum_{\mathbb{N}}n_k\log p_k\geq 0\}$ where $p_k$ is the $k$-th prime. We call $\tau(\mathbf{n})=\sum_{\mathbb{N}}n_k\log p_k$ the {\it ordinal number} of $\mathbf{n}$ and use $|\mathbf{n}|$ to denote the element of $\mathbb{Z}^\infty$ such that $\tau(|\mathbf{n}|)=|\tau(\mathbf{n})|$, i.e., $|\mathbf{n}|:=\mathbf{n}$ if $\tau(\mathbf{n})>0$ and $|\mathbf{n}|:=-\mathbf{n}$ if $\tau(\mathbf{n})<0$.
\end{example}

\begin{definition}\label{d3.3}
Suppose $G$ is a topological group with ordered dual $\Gamma$ where the order on $\Gamma$ is corresponding to the fixed semigroup $\Gamma'$ as in Definition \ref{d3.1}, we call $f\in L^1(G)$ an {\it analytic function} if $\supp(\hat{f})\subset \Gamma'$.
\end{definition}
\begin{definition}
The {\it Hardy space} $H^p(G)$ for $p\in [1,\infty]$ on $G$ is defined as the set of analytic $L^p(G)$-functions and $A(G)$ is defined as the set of all the continuous analytic functions on $G$. Naturally, the norms on these two spaces are the $L^p$-norm and the uniform norm, respectively.
\end{definition}
Throughout the paper, we use $\mathcal{P}$ to denote the set of analytic trigonometric polynomials, $\mathcal{P}_\mathbf{0}$ to denote trigonometric polynomials whose Fourier transform is supported on $\{\mathbf{n}\in \Gamma: \mathbf{n}>\mathbf{0}\}$, $[\mathcal{P}f]$ to denote the closed subspace of $H^2(G)$ generated by $\{pf: p\in \mathcal{P}\}$ and $[\mathcal{P}_\mathbf{0} f]$ in a similar way.

\begin{definition}\label{d3.5}
We call $f\in H^2(G)$ an outer function if $$e^{\int_G \log |f|dm_G}=\left|\int_G f dm_G\right|.$$
\end{definition}

In order to facilitate the exposition, we introduce some fundamental theorems of Fourier analysis on groups with ordered duals \cite[Theorem 8.3.2\& 8.4.3\& 8.5.2]{rud}.

\begin{theorem}[Szeg\"{o}]\label{t1.9}
Suppose $G$ is a compact group with ordered dual, $\mu$ is a finite Borel measure on $G$ with $\mu\geq 0$, and $d\mu=wdm_G+d\mu_s$ is the Lebesgue decomposition of $\mu$ with respect to $m_G$. Then
$$e^{\int_G \log w(x)dm_G}=\inf_{Q\in \mathcal{P}_\mathbf{0}}\int_G|1+Q(x)|^2d\mu(x).$$
\end{theorem}

\begin{theorem}\label{t1.10}
Suppose $G$ is a topological group with ordered dual, $w\in L^1(G)$ and $w\geq 0$. Then $w=|f|^2$ for $f\in H^2(G)$ with $\hat{f}(\mathbf{0})\neq 0$ if and only if $\int_{G}\log w(x)dx>-\infty$.
\end{theorem}

\begin{theorem}[Beurling, Helson and Lowdenslager]\label{t1.11}
Suppose $f\in H^2(G)$ and $$\int_{G}\log f(x)dx>-\infty,$$ Then $f=f_0f_1$, where $f_0$ is an inner function, i.e., $f_0\in H^2(G)$ and $|f_0|=1$ and $f_1$ is an outer function. The factorization is unique, except for multiplication by constants of absolute value 1.

Furthermore, $[\mathcal{P}f]=H^2(G)$ if and only if $f_1$ is an outer function.
\end{theorem}

It is well known that the dual of $\mathbb{T}^\omega$,  the Cartesian product of countable copies of $\mathbb{T}$, is $\mathbb{Z}^\infty$. Define $\mathbb{Z}_+^\infty$ be the smallest semigroup generated by countable copies of $\mathbb{N}$, and the positive semigroup $\{\mathbf{n}\in \mathbb{Z}^\infty: \tau(\mathbf{n})\geq 0\}$ of $\mathbb{Z}^\infty$ by $\mathbb{Z}^\infty_{\tau^+}$. It is clear that $\mathbb{Z}^\infty_+\in \mathbb{Z}^\infty_{\tau^+}$, and this allows us to view $H^p_+(\mathbb{T}^\omega)$ as a subspace of $H_\tau^p(\mathbb{T}^\omega)$.

\begin{definition}
Define $H^2_+(\mathbb{T}^\omega):=\{f\in L^2(\mathbb{T}^\omega): \supp(\hat{f})\subset \mathbb{Z}^\infty_+.\}$, $\mathcal{P}^+$ be all the trigonometric polynomials in $H^2_+(\mathbb{T}^\omega)$ and $\mathcal{P}_\mathbf{0}^+:=\mathcal{P}_\mathbf{0}\cap \mathcal{P}^+$, where $\mathcal{P}_\mathbf{0}$ is defined by the order from example \ref{x1.1}.
\end{definition}

For convenience, we will define some notations here that will be used in each subsequent section. In this paper, the bold-faced letter $\mathbf{a}$ always represents a number sequence, and its $k$-th coordinate is denoted by $a_k$, In particular, $\mathbf{1}$ represents a sequence where every coordinate is $1$ and $\mathbf{0}$ is defined similarly; $\mathbf{p}$ represents the sequence of prime numbers $\{2,3,5,...\}$. We use $\mathbf{n}$ to denote an element of $\mathbb{Z}^\infty$, $\gamma$ to denote a character of $\mathbb{T}^\omega$ and $\gamma_\mathbf{n}$ to denote the character according to $\mathbf{n}\in \mathbb{Z}^\infty$. For number sequences $\mathbf{a}$, $\mathbf{b}$, we use $\mathbf{a}\mathbf{b}$ to denote the sequence $\{a_kb_k\}$, and $\mathbf{a}^{\mathbf{b}}$ to denote $\prod a_k^{b_k}$. 

\section{A new proof of Szeg\"{o}'s theorem in infinitely many variables}
 The following known theorem was first proved by Szeg\"{o} on the unit disk \cite{s}:
\begin{theorem}
Suppose $\mu$ is a finite complex Borel measure on $\mathbb{T}$. $d\mu=wdm_{\mathbb{T}}+d\mu_s$ is the Lebesgue decomposition of $\mu$. Then
$$e^{\int_G\log w(x)dm_\mathbb{T}}=\inf_{Q\in \mathcal{P}_\mathbf{{0}}}\int_\mathbb{T} |1+Q(x)|^2 d\mu(x)$$
\end{theorem}
This theorem, like other theorems related to analytic functions, has different generalizations for groups with ordered dual and for plane domain. The case of groups has already been stated in Theorem \ref{t1.9}. As for the case of a plane domain, well-established results already exist. In 1991, the 2-dimensional case of this theorem was extended by Nakazi \cite{n}:
\begin{theorem}
Suppose $w\in L^1(\mathbb{T}^2), w\geq 0$ and $\log w\in L^1(\mathbb{T}^2)$. Then 
\begin{align}\label{e1.3}
\inf_{p\in \mathcal{P}_\mathbf{0}^+} \int_{\mathbb{T}^2}|1-p|^2 wdm_{\mathbb{T}^2}\geq e^{\int_{\mathbb{T}^2} \log wdm_{\mathbb{T}^2}}.
\end{align}
The equality holds if and only if there exists a cyclic vector $f$ for $H^2_+(\mathbb{T}^2)$ with $w=|f|^2$. i.e., $f\in H^2_+(\mathbb{T}^\omega)$ and $[\mathcal{P}^+f]=H^2_+(\mathbb{T}^2)$.
\end{theorem}

In 2022, Guo, Ni and Zhou extend this theorem to $H^q_+(\mathbb{T}^\omega),\ 1<q<\infty$ \cite[Theorem 3.1]{gnq}:

\begin{theorem}\label{t1.8}
Suppose $w\in L^1(\mathbb{T}^\omega), w\geq 0$ and $\log w\in L^1(\mathbb{T}^\omega)$. Then 
\begin{align}\label{e4.1}
\inf_{p\in \mathcal{P}_\mathbf{0}^+} \int_{\mathbb{T}^\omega}|1-p|^q wdm_{\mathbb{T}^\omega}\geq e^{\int_{\mathbb{T}^\omega} \log wdm_{\mathbb{T}^\omega}}
\end{align}
The equality holds if and only if there exists a cyclic vector $f$ for $H^q_+(\mathbb{T}^\omega)$ with $w=|f|^q$. 
\end{theorem}

In order to demonstrate the usefulness of the order $\tau$, the present section is devoted to showing a brief proof of Theorem \ref{t1.8} in the case of $q=2$. The first part of this theorem, namely inequality \eqref{e4.1}, can be viewed as a direct corollary of Theorem \ref{t1.9} since $\mathcal{P}^+_\mathbf{0}\subset \mathcal{P}_\mathbf{0}.$

The following proposition is quite useful in our proofs.

\begin{proposition}\label{p1.1}
Suppose $G$ is a locally compact abelian group and $f, g\in L^2(G)$, then $\supp(\hat{fg})\in \supp(\hat{f})+\supp(\hat{g})$.
\end{proposition}
\begin{proof}
By the Plancherel theorem,  there are sequences of trigonometric polynomials $\{f_k\}, \{g_k\}$ converging to $f$, $g$ in $L^2(G)$ respectively, with $\supp (f_k)\subset \supp(f)$ and $\supp(g_k)\subset \supp(g)$, thus $\supp(\hat{f_kg_k})\subset \supp(\hat{f}\ast\hat{g})\subset \supp(\hat{f})+\supp(\hat{g})$. Since $f_kg_k\to fg$ in $L^1(G)$, it follows that $\hat{f_kg_k}\to \hat{fg}$ uniformly on the dual group of $G$, thus $\supp(\hat{fg})\subset \supp(\hat{f})+\supp(\hat{g})$.
\end{proof}

The following proof borrowed some idea from that of Nakazi \cite{n}.

\begin{proof}[Proof of the second part of Theorem \ref{t1.8}]
We start with the case that $f$ is a cyclic vector for $H^2_+(\mathbb{T}^\omega)$. Indeed, $f$ is also a cyclic vector for $H_{\tau}^2(\mathbb{T}^\omega)$ and $\int f\neq 0$ since the unit function belongs to $[\mathcal{P}^+f]$ and $\mathcal{P}^+\subset \mathcal{P}$. By Theorem \ref{t1.11}, $f$ is an outer function, so we obtain
\begin{align}\label{e3.1}
|\hat{f}(\mathbf{0})|^2=\left|\int_{\mathbb{T}^\omega} f dm_{\mathbb{T}^\omega}\right|^2=e^{\int_{\mathbb{T}^\omega} \log |f|^2dm_{\mathbb{T}^\omega}}.
\end{align}

By Proposition \ref{p1.1}, $\supp(\hat{fp})\subset \supp(\hat{p})+\supp(\hat{f})\subset Z^\infty_{\tau^+}\setminus \{\mathbf{0}\}$ for every $p\in \mathcal{P}_\mathbf{0}^+$, together with Plancherel's theorem we obtain that
\begin{align}\label{e1.1}
\inf_{p\in \mathcal{P}_\mathbf{0}^+} \int_{\mathbb{T}^\omega}|f-fp|^2dm_{\mathbb{T}^\omega}  \geq |\hat{f}(\mathbf{0})|^2.
\end{align}
Since $f$ is a cyclic vector for $H^2_+(\mathbb{T}^\omega)$, there are trigonometric polynomials $\{p_k\}\subset H^2_+(\mathbb{T}^\omega)$ such that 
\begin{align}\label{e4.6}
p_kf\to f-\hat{f}(\mathbf{0})
\end{align}
with respect to $L^2$-norm as $k\to \infty$. Combining this with $\hat{f}(\mathbf{0})\neq 0$ we have $\hat{p_k}(\mathbf{0})\to 0$ which allows us to replace $p_k$ by $p_k-\hat{p_k}(\mathbf{0})$ without affecting \eqref{e4.6}, thus it follows that
\begin{align}\label{e1.2}
\inf_{p\in \mathcal{P}_\mathbf{0}^+} \int_{\mathbb{T}^\omega}|f-fp|^2dm_{\mathbb{T}^\omega}\le \lim_{k\to \infty} \int_{\mathbb{T}^\omega}|f-f(p_k-\hat{p_k}(0))|^2dm_{\mathbb{T}^\omega}=|\hat{f}(\mathbf{0})|^2.
\end{align}
From \eqref{e3.1}, \eqref{e1.1} and \eqref{e1.2}, we deduce that equality \eqref{e4.1} holds.

Conversely, suppose $w$ satisfies $$\inf_{p\in \mathcal{P}_\mathbf{0}^+} \int_{\mathbb{T}^\omega}|1-p|^2 w dm_{\mathbb{T}^\omega}=e^{\int_{\mathbb{T}^\omega} \log wdm_{\mathbb{T}^\omega}}.$$ By Theorem \ref{t1.10} and Theorem \ref{t1.11}, we have that $w=|f|^2$ for some outer function $f\in H_{\tau}^2(\mathbb{T}^\omega)$. That is,  $$\inf_{p\in \mathcal{P}_\mathbf{0}^+} \int_{\mathbb{T}^\omega}|f-fp|^2 dm_{\mathbb{T}^\omega}=e^{2\int_{\mathbb{T}^\omega} \log |f| dm_{\mathbb{T}^\omega}}=|\hat{f}(\mathbf{0})|^2.$$ Write $f=f_0+f_1$ where $f_0$ lies in $[\mathcal{P}^+f]\ominus[\mathcal{P}_\mathbf{0}^+f]$ and $f_1\in [\mathcal{P}_\mathbf{0}^+f]$. By Proposition \ref{p1.1}, $\supp(\hat{pf})\subset \supp(\hat{p})+\supp(\hat{f})\subset Z^\infty_{\tau^+}\setminus \{\mathbf{0}\}$ for every $p\in \mathcal{P}_\mathbf{0}^+$,
we have $$|\hat{f_0}(\mathbf{0})|^2=|\hat{f}(\mathbf{0})|^2=\inf_{p\in \mathcal{P}_\mathbf{0}^+} \int_{\mathbb{T}^\omega}|f-fp|^2 dm_{\mathbb{T}^\omega}=\|f_0\|^2,$$ By Plancherel's theorem, $f_0$ is a constant and thus $H^2_+(\mathbb{T}^\omega)\subset[\mathcal{P}^+f]$.

It remains to show that $f\in H^2_+(\mathbb{T}^\omega)$. Since $H^2_+(\mathbb{T}^\omega)\subset[\mathcal{P}^+f]$ there exist trigonometric polynomials $\{p_k\}\subset \mathcal{P}^+$ such that $p_kf\to \hat{f}(\mathbf{0})\bold{1}$ with respect to $L^2$-norm. Since $\supp(\hat{p_k})\subset \mathbb{Z}^\infty_+\subset \mathbb{Z}^\infty_{\tau^+}$ and $\supp(\hat{f})\subset \mathbb{Z}^\infty_{\tau^+}$, by Proposition \ref{p1.1} it follows $\hat{p_kf}(\mathbf{0})=\hat{p_k}(\mathbf{0})\hat{f}(\mathbf{0})$ and $\hat{p_k}(\mathbf{0})\to 1$. Without loss of generality, we may assume that $\hat{p_k}(\mathbf{0})=1$.

Suppose $\supp(\hat{f})\setminus \mathbb{Z}^\infty_+$ is non-empty, select $\mathbf{n}_0$ in $\supp(\hat{f})\setminus \mathbb{Z}^\infty_+$ such that $$\tau(\mathbf{n}_0)<\log 2+\inf_{\mathbf{n}\in \supp(\hat{f})\setminus \mathbb{Z}^\infty_+}\tau(\mathbf{n})$$ and let $f_{\mathbf{n}_0}$ be the $\mathbf{n}_0$-term in Fourier series of $f$. We have
\begin{align*}
p_kf=\hat{f}(\mathbf{0})+f_{\mathbf{n}_0}+(p_k-1)f+(f-\hat{f}(\mathbf{0})-f_{\mathbf{n}_0}).
\end{align*}
\qquad Since the ordinal number of each term in the Fourier series of $p_k-1$ is greater than or equal to $\log 2$ and $\mathbb{Z}^\infty_+ +\mathbb{Z}^\infty_+\subset \mathbb{Z}^\infty_+$, we have $\supp (\hat{(p_k-1)f)}\cap (\{\mathbf{0}\}\cup \{\mathbf{n}_0\})=\emptyset$, by Plancherel's theorem we have $\|p_kf\|_2^2\geq |\hat{f}(\mathbf{0})|^2+|\hat{f}(\mathbf{n}_0)|^2$ which lead a contradiction to $p_kf\to \hat{f}(\mathbf{0})\bold{1}$.
\end{proof}

\begin{remark}
Guo et al. demonstrated the proof of this theorem by employing their generalized theory of Smirnov function classes \cite[Proposition 2.11]{gnq}. Such an approach holds greater advantages for the general $H^p$ case. In contrast, the Fourier-analytic method we adopted relies on orthogonal decomposition in Hilbert spaces and is only readily applicable in the $H^2$ setting.
\end{remark}

\section{Poisson measures on $\mathbb{T}^\omega$}

\quad Analytic functions on the closed unit disk have an important characteristic that they can be equivalently defined from two perspectives: the boundary and the interior of the disk. The reason why these functions on the circle can be extended into the interior of the disk is that the maximal ideal space of the disk algebra, i.e., the algebra of analytic continuous functions on $\mathbb{T}$,
is richer than $\mathbb{T}$. This allows the interior of the disk to be viewed as a subset of the maximal ideal space of the disk algebra.

A natural question is whether the maximal ideal space of $A_\tau(\mathbb{T}^\omega)$ serves as the counterpart of the one-dimensional unit disk in the infinite-dimensional case? Specifically, can this maximal ideal space be naturally associated to a certain set in $\mathbb{C}^\omega$ (Cartesian product of countable copies of $\mathbb{C}$) whose topological boundary is $\mathbb{T}^\omega$. To answer this question, we need to find all the complex homomorphisms of $A_\tau(\mathbb{T}^\omega)$ and the corresponding representing measures. In the classical case, such a measure is induced by a Poisson kernel.

Our investigation of this maximal ideal space follows conventional approaches. In the case of topological groups (as opposed to semigroups), the Wiener algebra on $\Gamma$, which consists of Fourier transform of the $L^1(G)$-functions, is dense in $C_0(\Gamma)$. Moreover, its topology is stronger than $L^1(G)$ because the Fourier transform is continuous, which implies $L^1(G)$ and it's Wiener algebra sharing identical maximal ideal spaces. In our discussion, $L^1(G)$ and $C_0(\Gamma)$ are replaced by $l_1(\mathbb{Z}^\infty_{\gamma_{\tau^+}})$ and $A_\tau(\mathbb{T}^\omega)$, respectively.

From now on we use $\mathbf{r}$ to denote a sequence of positive real numbers $\mathbf{r}=\{r_1,r_2,...,r_k,...\}$ with $r_1\le 1$.

%
%
%
%
%
%
%
%
%
%
%

\begin{lemma}\label{l1.1}
For $\boldsymbol{\lambda}\in \mathbb{T}^\omega$, $(\mathbf{r}\boldsymbol{\lambda})^\mathbf{n}$ is bounded for all $\mathbf{n}\in \mathbb{Z}^\infty_{\tau^+}$ if and only if for all $k>1$, $\log_{p_k}r_k=\log_2r_1$, where $p_k$ is the $k$-th prime number.
\end{lemma}

\begin{proof}

Since $|(\mathbf{r}\boldsymbol{\lambda})^\mathbf{n}|=e^{\sum n_k\log r_k}$ and $\mathbb{Z}^\infty_{\tau^+}$ is closed under positive scalar multiplication, we just need to prove that $\sum n_k\log r_k\le 0$ holds for all $\mathbf{n}\in \mathbb{Z}^\infty_{\tau^+}$ if and only if $\log r_k=   \frac{\log\frac{1}{r_1}}{\log 2}\log \frac{1}{p_k}$. Recall that $\mathbf{n}\in \mathbb{Z}^\infty_{\tau^+}$ means that $\sum n_k \log p_k\geq 0$. It is clear that if $\log r_k= \frac{\log \frac{1}{r_1}}{\log 2}\log \frac{1}{p_k}$ then $(\mathbf{r}\boldsymbol{\lambda})^\mathbf{n}$ is bounded.

Conversely, we suppose that $\sum n_k\log r_k\le 0$ for all $\mathbf{n}\in \mathbb{Z}^\infty_{\tau^+}$. By taking $\mathbf{n}$ as such an element in $\mathbb{Z}^\infty$ that one of its coordinate is 1 and the others are 0, one can easily see that $\log \frac{1}{r_k}$ is nonnegative for all $k\in \mathbb{N}$. For arbitrary $k$, choose $\{a_m\}, \{b_m\}\subset \mathbb{Z}$ such that $a_m\log 2+b_m\log p_k\to 0^+$ with $b_m\geq 0$ for all $m\in \mathbb{N}$. By our assumption, $a_m\log\frac{1}{r_1}+b_m\log\frac{1}{r_k}\geq 0$. From this we deduce that $\frac{\log\frac{1}{r_k}}{\log\frac{1}{r_1}}\geq \frac{\log p_k}{\log 2}$. The reverse inequality can be obtained by assuming $b_m\le 0$ for all $m\in \mathbb{N}$, which completes this proof.
\end{proof}

\begin{corollary}\label{c1.1}
The maximal ideal space $\Delta$ of $l_1(\mathbb{Z}^\infty_{\tau^+})$ coincides with the set $$\mathfrak{D}:=\left\{\mathbf{r}\boldsymbol{\lambda}: \boldsymbol{\lambda}\in \mathbb{T}^\omega, \log_{p_k}r_k=\log_2r_1\right\},$$
namely, each complex homomorphism of $l_1(\mathbb{Z}^\infty_{\tau^+})$ takes the form:
$$f\to \sum_{\mathbf{n}\in \mathbb{Z}^\infty_{\tau^+}}f(\mathbf{n})(\mathbf{r}\boldsymbol{\lambda})^\mathbf{n}.$$
\end{corollary}

\begin{proof}
This proof is almost identical to \cite[Theorem 1.2.2]{rud}.

For a given $(\mathbf{r}\boldsymbol{\lambda})^\mathbf{n}\in \mathfrak{D}$ and $f,g\in l_1(\mathbb{Z}^\infty_{\tau^+})$, we have
$$\sum_{\mathbf{n}\in\mathbb{Z}^\infty_{\tau^+}}f\ast g(\mathbf{n})(\mathbf{r}\boldsymbol{\lambda})^\mathbf{n}=(\sum_{\mathbf{n}\in\mathbb{Z}^\infty_{\tau^+}}f(\mathbf{n})(\mathbf{r}\boldsymbol{\lambda})^\mathbf{n})(\sum_{\mathbf{n}\in\mathbb{Z}^\infty_{\tau^+}}g(\mathbf{n})(\mathbf{r}\boldsymbol{\lambda})^\mathbf{n}).$$
Thus the map $f\mapsto \sum_{\mathbf{n}\in\mathbb{Z}^\infty_{\tau^+}}f(\mathbf{n})(\mathbf{r}\boldsymbol{\lambda})^\mathbf{n}$ induced by $(\mathbf{r}\boldsymbol{\lambda})^\mathbf{n}$ is a complex homomorphism.

Conversely, for a given complex homomorphism $\chi\in \Delta$, by the Riesz representation theorem there is a $l^\infty(\mathbb{Z}^\infty_{\tau^+})$-function $\phi_\chi$ such that $\chi(f)=\sum_{\mathbf{n}\in \mathbb{Z}^\infty_{\tau^+}} f(\mathbf{n})\phi_\chi(\mathbf{n})$. Since $\chi(f\ast g)=\chi(f)\chi(g)$ we have $\phi_\chi(\mathbf{n}+\mathbf{m})=\phi_\chi(\mathbf{n})\phi_\chi(\mathbf{m})$, and thus $\phi_\chi(\mathbf{n})$ takes the form $\phi_\chi(\mathbf{n})=(\mathbf{r}\boldsymbol{\lambda})^\mathbf{n}$ where $\mathbf{r}, \boldsymbol{\lambda}$ are corresponding to $\chi$. Indeed, if we use $\{\delta_k\}$ to denote the sequence whose $k$-th coordinate is $1$ and others are $0$, then $r_k\boldsymbol{\lambda}_k=\phi_\chi(\delta_k)$.

From $\phi_\chi\in l^\infty(\mathbb{Z}^\infty_{\tau^+})$ and Lemma \ref{l1.1}, we deduce that $\mathbf{r}\boldsymbol{\lambda}\in \mathfrak{D}$. 
\end{proof}

\begin{remark}\label{r1.1}
For convenience, we use $\mathbf{p}^{-\sigma}$ to denote the sequence $\{2^{-\sigma},3^{-\sigma},5^{-\sigma},...\}$ for $\sigma\in [0,\infty)$, and $\mathbf{0}$ for $\sigma=\infty$. Then the definition of $\mathfrak{D}$ can be rewritten as
$$\mathfrak{D}:=\{\mathbf{p}^{-\sigma}\boldsymbol{\lambda}: \boldsymbol{\lambda}\in \mathbb{T}^\omega, \sigma\in [0,\infty]\}.$$
From this one can easily see that there is a bijection from $\mathfrak{D}$ to  $([0,\infty)\times 
\mathbb{T}^\omega)\cup\{\boldsymbol{\infty}\}:
\mathbf{r}\boldsymbol{\lambda}\mapsto (\sigma, \boldsymbol{\lambda}),\ \mathbf{r}=\mathbf{p}^{-\sigma}\neq \mathbf{0}$ and $\mathbf{0}\mapsto \boldsymbol{\infty}$, which allows us to denote each element $\mathbf{r}\boldsymbol{\lambda}$ in $\mathfrak{D}$
by $\sigma \cdot \boldsymbol{\lambda}$.
Furthermore, following a similar approach to the Alexandrov topology, we may equip the set $([0,\infty)\times 
\mathbb{T}^\omega)\cup\{\boldsymbol{\infty}\}$ (or equivalently, the set $\mathfrak{D}$) with a new topology to make it a compactification of $[0,\infty)\times \mathbb{T}^\omega$, whose members are those following sets:
\begin{itemize}
    \item 1.  the open sets in $[0,\infty)\times \mathbb{T}^\omega$ with respect to product topology.
    \item 2.  the sets containing $\boldsymbol{\infty}$ and whose complement has the form $K\times \mathbb{T}^\omega$, where $K$ is an arbitrary compact set in $[0,\infty)$.
\end{itemize}

\end{remark}
It is easy to check that $\hat{l_1}(\mathbb{Z}^\infty_{\tau^+}):=\{\hat{f}: f\in l_1(\mathbb{Z}^\infty_{\tau^+})\}$
 is a subalgebra of $A_\tau(\mathbb{T}^\omega)$, where $\hat{f}$ is the Fourier transform of $f$. We ask if $\mathfrak{D}$ can serve as a subset of the maximal ideal space of $A_\tau(\mathbb{T}^\omega)$. According to the Hahn-Banach theorem, every complex homomorphism in $\mathfrak{D}$ can be represented by a complex measure on $\mathbb{T}^\omega$. In fact, we can further require such measures to be positive, which will provide significant convenience for subsequent discussions.

\begin{lemma}\label{l1.2}
Suppose $A$ is a finite subset of $\mathbb{Z}^\infty$, and $\{c_\mathbf{n}\}_{\mathbf{n}\in \mathbb{Z}^\infty}$ is a sequence of complex numbers such that $c_\mathbf{n}=0, \forall \mathbf{n}\notin A$, it follows
\begin{align}\label{e1.5}
\sum_{\mathbf{n},\mathbf{m}\in \mathbb{Z}^\infty} c_\mathbf{n}\bar{c_\mathbf{m}}\mathbf{r}^{|\mathbf{n}-\mathbf{m}|}\geq 0.
\end{align}
where $\mathbf{n}\in \mathbb{Z}^\infty$, $\mathbf{r}\in \mathfrak{D}$.
\end{lemma}

\begin{proof}
Without loss of generality,  we assume that the finite set $A\subset \mathbb{Z}^\infty_+$. Indeed, for arbitrary $\mathbf{k}\in \mathbb{Z}^\infty$, if we replace $\{c_\mathbf{n}\}_{\mathbf{n}\in \mathbb{Z}^\infty}$ by $\{d_\mathbf{n}\}_{\mathbf{n}\in \mathbb{Z}^\infty}$ with  $d_\mathbf{n}=c_{\mathbf{n}+\mathbf{k}}$, the sum in \eqref{e1.5} will not change.

We use $\mathbf{n}_k$ to denote the k-th element of $A$ with respect to order $\tau$ and let $N+1:=card(A)$. It is sufficient to prove that the following matrix
\begin{align}\label{e1.6}
\left[
\begin{matrix}
&1 &\mathbf{r}^{\mathbf{n}_1-\mathbf{n}_0} &\mathbf{r}^{\mathbf{n}_2-\mathbf{n}_0}&...&\mathbf{r}^{\mathbf{n}_N-\mathbf{n}_0}\\
&\mathbf{r}^{\mathbf{n}_1-\mathbf{n}_0} &1 &\mathbf{r}^{\mathbf{n}_2-\mathbf{n}_1} &...&\mathbf{r}^{\mathbf{n}_N-\mathbf{n}_1}\\
&\mathbf{r}^{\mathbf{n}_2-\mathbf{n}_0} &\mathbf{r}^{\mathbf{n}_2-\mathbf{n}_1} &1&...&\mathbf{r}^{\mathbf{n}_N-\mathbf{n}_2}\\
&\vdots &\vdots &\vdots & &\vdots\\
&\mathbf{r}^{\mathbf{n}_N-\mathbf{n}_0} &\mathbf{r}^{\mathbf{n}_N-\mathbf{n}_1} &\mathbf{r}^{\mathbf{n}_N-\mathbf{n}_2} &\ldots &1
\end{matrix}
\right]
\end{align}
is positive for each $N\in \mathbb{N}$:

\begin{align}\label{e1.7}\nonumber
&\left|
\begin{matrix}
&1 &\mathbf{r}^{\mathbf{n}_1-\mathbf{n}_0} &\mathbf{r}^{\mathbf{n}_2-\mathbf{n}_0}&...&\mathbf{r}^{\mathbf{n}_N-\mathbf{n}_0}\\
&\mathbf{r}^{\mathbf{n}_1-\mathbf{n}_0} &1 &\mathbf{r}^{\mathbf{n}_2-\mathbf{n}_1} &...&\mathbf{r}^{\mathbf{n}_N-\mathbf{n}_1}\\
&\mathbf{r}^{\mathbf{n}_2-\mathbf{n}_0} &\mathbf{r}^{\mathbf{n}_2-\mathbf{n}_1} &1&...&\mathbf{r}^{\mathbf{n}_N-\mathbf{n}_2}\\
&\vdots &\vdots &\vdots & &\vdots\\
&\mathbf{r}^{\mathbf{n}_N-\mathbf{n}_0} &\mathbf{r}^{\mathbf{n}_N-\mathbf{n}_1} &\mathbf{r}^{\mathbf{n}_N-\mathbf{n}_2} &\ldots &1
\end{matrix}
\right|\\ \nonumber
=&
\left|
\begin{matrix}
&1 \ &\mathbf{r}^{\mathbf{n}_1-\mathbf{n}_0} \ &\mathbf{r}^{\mathbf{n}_2-\mathbf{n}_0}\ &...\ &\mathbf{r}^{\mathbf{n}_N-\mathbf{n}_0}\\
&\mathbf{r}^{\mathbf{n}_1-\mathbf{n}_0}-\mathbf{r}^{\mathbf{n}_0-\mathbf{n}_1} \ &0 \ &0 \ &...\ &0\\
&\mathbf{r}^{\mathbf{n}_2-\mathbf{n}_0}-\mathbf{r}^{\mathbf{n}_0-\mathbf{n}_2} \ &\mathbf{r}^{\mathbf{n}_2-\mathbf{n}_1}-\mathbf{r}^{\mathbf{n}_1-\mathbf{n}_2} \ &0\ &...\ &0\\
&\vdots \ &\vdots \ &\vdots \ & \ &\vdots\\
&\mathbf{r}^{\mathbf{n}_N-\mathbf{n}_0}-\mathbf{r}^{\mathbf{n}_0-\mathbf{n}_N} \ &\mathbf{r}^{\mathbf{n}_N-\mathbf{n}_1}-\mathbf{r}^{\mathbf{n}_1-\mathbf{n}_N} \ &\mathbf{r}^{\mathbf{n}_N-\mathbf{n}_2}-\mathbf{r}^{\mathbf{n}_2-\mathbf{n}_N} \ &\ldots \ &0
\end{matrix}
\right|\\ \nonumber
=&
\left|
\begin{matrix}
&1 \ &\mathbf{r}^{\mathbf{n}_1-\mathbf{n}_0} \ &\mathbf{r}^{\mathbf{n}_2-\mathbf{n}_0}\ &...\ &\mathbf{r}^{\mathbf{n}_N-\mathbf{n}_0}\\
&\mathbf{r}^{\mathbf{n}_1-\mathbf{n}_0}-\mathbf{r}^{\mathbf{n}_0-\mathbf{n}_1} \ &0 \ &0 \ &...\ &0\\
&0 \ &\mathbf{r}^{\mathbf{n}_2-\mathbf{n}_1}-\mathbf{r}^{\mathbf{n}_1-\mathbf{n}_2} \ &0\ &...\ &0\\
&\vdots \ &\vdots \ &\vdots \ & \ &\vdots\\
&0 \ &0 \ &0 \ &\ldots \ &0
\end{matrix}
\right|\\ \nonumber
=&
\left|
\begin{matrix}
&0 \ &0 \ &0\ &...\ &\mathbf{r}^{\mathbf{n}_N-\mathbf{n}_0}\\
&\mathbf{r}^{\mathbf{n}_1-\mathbf{n}_0}-\mathbf{r}^{\mathbf{n}_0-\mathbf{n}_1} \ &0 \ &0 \ &...\ &0\\
&0 \ &\mathbf{r}^{\mathbf{n}_2-\mathbf{n}_1}-\mathbf{r}^{\mathbf{n}_1-\mathbf{n}_2} \ &0\ &...\ &0\\
&\vdots \ &\vdots \ &\vdots \ & \ &\vdots\\
&0 \ &0 \ &0 \ &\ldots \ &0
\end{matrix}
\right|\\
=&\mathbf{r}^{\mathbf{n}_N-\mathbf{n}_0}\prod_{k=0}^N (\mathbf{r}^{\mathbf{n}_k-\mathbf{n}_{k+1}}-\mathbf{r}^{\mathbf{n}_{k+1}-\mathbf{n}_k})
\end{align}
By our assumption that $\tau(\mathbf{n}_{k+1}-\mathbf{n}_k)\geq 0$, we know that $\mathbf{r}^{\mathbf{n}_{k+1}-\mathbf{n}_k}=e^{\sigma\tau(\mathbf{n}_{k+1}-\mathbf{n}_k)}\geq 1$, where $\sigma$ is as in Remark \ref{r1.1}. Thus \eqref{e1.7} is positive for $\mathbf{r}\in \mathfrak{D}$, and this implies that the matrix \eqref{e1.6} is positive for $\mathbf{r}\in \mathfrak{D}$.
%
%
%
%
\end{proof}

\begin{definition}\label{d5.5}
By Bochner's theorem \cite[Theorem 1.4.3]{rud} and Lemma \ref{l1.2}, if $\mathbf{r}\in \mathfrak{D}$ then there is a unique positive Borel measure $P_\mathbf{r}$ on $\mathbb{T}^\omega$ with $\hat{P}(\mathbf{n})=\mathbf{r}^{|\mathbf{n}|}$ for every $\mathbf{n}\in \mathbb{Z}^\infty$. We call $P_\mathbf{r}$ the {\it Poisson measure} on $\mathbb{T}^\omega$ according to $\mathbf{r}$ and use $P_{\mathbf{r},\boldsymbol{\lambda}}$ to denote the measure such that $P_{\mathbf{r},\boldsymbol{\lambda}}(E)=P_\mathbf{r}(E\boldsymbol{\lambda}^{-1})$ for every Borel set $E$ on $\mathbb{T}^\omega$. 
\end{definition}

The Poisson measure is of norm $1$ since $\hat{P}_\mathbf{r}(0)=1$. Furthermore, $P_{\mathbf{r}}$ is singular to $m_{\mathbb{T}^\omega}$ for each $\mathbf{r}\in \mathfrak{D}\setminus\{\mathbf{0}\}$. For if it is false, we write $P_{\mathbf{r}}:=\sigma_\mathbf{{r}}+s_{\mathbf{r}}$ where $\sigma_\mathbf{{r}}$ and $s_{\mathbf{r}}$ are absolutely continuous part and singular part of $P_\mathbf{r}$, respectively. $\sigma_\mathbf{r}$ and $s_{\mathbf{r}}$ are both positive. By generalized Riemann-Lebesgue theorem, for each $\epsilon$ there is a finite subset $F$ of $\mathbb{Z}^\infty$ such that $\sup_{\mathbf{n}\in \mathbb{Z}^\infty\setminus F} \hat{\sigma_\mathbf{r}}(\mathbf{n})\le \epsilon$. This gives $\sup_{n\in \mathbb{Z}^\infty}\hat{s_\mathbf{r}}(\mathbf{n})=1$. However, if $\hat{\sigma_{r}}(\mathbf{0})>0$, then $\hat{s_{\mathbf{r
}}}(\mathbf{0})<1$, contradicting to the positivity of $s_{\mathbf{r}}$.

There is a useful result for $A(G)$ on arbitrary compact abelian group $G$ \cite[8.7.3]{rud}. For convenience, we will present its proof here.
\begin{proposition}\label{p5.6}
The set of analytic trigonometric polynomials is dense in $A(G)$ for arbitrary compact abelian group $G$ with ordered dual.
\end{proposition}
\begin{proof}
By the Hahn-Banach theorem and the Riesz representation theorem, each linear continuous functional on $A(G)$ can be extended to a functional on $C(G)$ which coincides with a complex Borel measure $\mu$ on $G$. It is sufficient to prove each Borel measure $\mu$ on $G$ that annihilates every analytic trigonometric polynomial also annihilates $A(G)$. Suppose that $\hat{\mu}(\mathbf{n})=0$ for each $\mathbf{n}\le0$, then for each $g\in A(G)$, we have $\int_G g d\mu=\tilde{g}\ast \mu(\mathbf{0})$, where $\tilde{g}(x):=g(x^{-1})$. By the compactness of $G$, $\tilde{g}\ast \mu$ is continuous. Combining this with $\hat{\tilde{g}}\hat{\mu}=0$, we conclude that $\tilde{g}\ast \mu(\mathbf{0})=0$.
\end{proof}

The following theorems show that the Poisson measure on $\mathbb{T}^\omega$ is a natural generalization of that in one-variable setting.

\begin{theorem}\label{t1.5}
$P_{\mathbf{r},\boldsymbol{\lambda}}$ induces a complex homomorphism on the disk algebra $A_\tau(\mathbb{T}^\omega)$. 
\end{theorem}

\begin{proof}
Suppose $f,g$ are trigonometric polynomials on $\mathbb{T}^\omega$ with $\supp(f)\subset \mathbb{Z}^\infty_{\tau^+}, \supp(g)\subset \mathbb{Z}^\infty_{\tau^+}$, then by Corollary \ref{c1.1} we have \begin{align*}
\int_{\mathbb{T}^\omega}fg \ dP_{\mathbf{r},\boldsymbol{\lambda}}&=\sum_{\mathbf{n}\in \mathbb{Z}^\infty} \hat{f}\ast \hat{g}(\mathbf{n})(\mathbf{r}\boldsymbol{\lambda})^\mathbf{n}=\left(\sum_{\mathbf{n}\in \mathbb{Z}^\infty}\hat{f}(\mathbf{n})(\mathbf{r}\boldsymbol{\lambda})^\mathbf{n}\right)\left(\sum_{\mathbf{n}\in \mathbb{Z}^\infty} \hat{g}(\mathbf{n})(\mathbf{r}\boldsymbol{\lambda})^\mathbf{n}\right)\\
&=\left(\int_{\mathbb{T}^\omega}f dP_{\mathbf{r},\boldsymbol{\lambda}}\right)\left(\int_{\mathbb{T}^\omega}g dP_{\mathbf{r},\boldsymbol{\lambda}}\right).
\end{align*}

By Proposition \ref{p5.6}, for $f, g\in A_\tau(\mathbb{T}^\omega)$, there are $\{f_k\}_{k=1}^\infty, \{g_k\}_{k=1}^\infty$, which are sequences of polynomials uniformly tending to $f, g$ respectively, then by Proposition \ref{p1.1} we have
\begin{align*}
\int_{\mathbb{T}^\omega}fgdP_{\mathbf{r},\boldsymbol{\lambda}}&=\lim_{k\to\infty}\int_{\mathbb{T}^\omega}f_kg_kdP_{\mathbf{r},\boldsymbol{\lambda}}\\
&=\lim_{k\to\infty}\int_{\mathbb{T}^\omega}f_kdP_{\mathbf{r},\boldsymbol{\lambda}}\int_{\mathbb{T}^\omega}g_kdP_{\mathbf{r},\boldsymbol{\lambda}}\\
&=
\int_{\mathbb{T}^\omega}fdP_{\mathbf{r},\boldsymbol{\lambda}}\int_{\mathbb{T}^\omega}gdP_{\mathbf{r},\boldsymbol{\lambda}}.
\end{align*}
This completes the proof.
\end{proof}

\begin{corollary}
The maximal ideal space $\tilde{\Delta}$ of $A_\tau(\mathbb{T}^\omega)$ coincides with $\mathfrak{D}$.
\end{corollary}

\begin{proof}
Theorem \ref{t1.5} already shows that $\mathfrak{D}\subset\tilde{\Delta}$.
Note that:

(1) The Fourier transform from $l_1(\mathbb{Z}^\infty_{\tau^+})$ to $A_\tau(\mathbb{T}^\omega)$ is continuous.

(2) The set of all analytic trigonometric polynomials is dense in $A_\tau(\mathbb{T}^\omega)$ (Proposition \ref{p5.6}). 

Thus each complex homomorphism $\chi$ of $l_1(\mathbb{Z}^\infty_{\tau^+})$ uniquely induces a complex homomorphism $\chi'$ of $A_\tau(\mathbb{T}^\omega)$ such that
$\chi'(\hat{f})=\chi(f)$ for all $f\in l_1(\mathbb{{Z}^\infty_{\tau^+}})$, where $\hat{f}$ is as in Remark \ref{r1.1}. This shows that $\tilde{\Delta}\subset \mathfrak{D}$ which completes the proof.
\end{proof}

In the setting of the unit disk, functions in the $L^p(\mathbb{T}), 1\le p<\infty$ can be approximated by their harmonic extensions from the interior of the disk. This approximation property admits a natural extension to the case we consider here.

\begin{definition}
For $f\in L_p(\mathbb{T}^\omega)$, $1\le p<\infty$, define $$f_{\mathbf{r}}(\boldsymbol{\lambda})=f\ast P_\mathbf{r}(\boldsymbol{\lambda})=\int_{\mathbb{T}^\omega} f dP_{\mathbf{r},\boldsymbol{\lambda}}=\int_{\mathbb{T}^\omega} f(\boldsymbol{\lambda}\boldsymbol{\beta}^{-1}) dP_{\mathbf{r}}(\boldsymbol{\beta}).$$
\end{definition}

\begin{proposition}
\label{p1.6}
For $f\in L^p(\mathbb{T}^\omega)$, $1\le p<\infty$, it follows that $\|f_\mathbf{r}\|_p\le \|f\|_p$.
\end{proposition}

\begin{proof}
By generalized Minkowski inequality, for $p\in[1,\infty)$, we have
\begin{align*}
\left(\int_{\mathbb{T}^\omega} |f_\mathbf{r}|^p dm_{\mathbb{T}^\omega}\right)^{\frac{1}{p}}&=\left(\int_{\mathbb{T}^\omega}\left|\int_{\mathbb{T}^\omega}f(\boldsymbol{\lambda}\boldsymbol{\beta}^{-1})dP_{\mathbf{r}}(\boldsymbol{\beta})\right|^p dm_{\mathbb{T}^\omega}(\boldsymbol{\lambda})\right)^{\frac{1}{p}}\\
&\le\int_{\mathbb{T}^\omega}\left(\int_{\mathbb{T}^\omega}|f|^p dm_{\mathbb{T}^\omega}\right)^{\frac{1}{p}}d P_\mathbf{r}\\
&=\|f\|_p.
\end{align*}
\end{proof}

\begin{corollary}\label{c5.11}
For $f\in H_{\tau}^2(\mathbb{T}^\omega)$ and $g\in A_\tau(\mathbb{T}^\omega)$, $(fg)_\mathbf{r}=f_\mathbf{r}g_\mathbf{r}$.
\end{corollary}
\begin{proof}
Let $\{f_n\}_{n=1}^\infty$ be a sequence of trigonometric polynomials such that $\|f_n-f\|_2\to 0$, we have
\begin{align*}
\|(fg)_\mathbf{r}-f_\mathbf{r}g_\mathbf{r}\|_2&=\|(fg)_\mathbf{r}-(f_ng)_\mathbf{r}+(f_ng)_\mathbf{r}-f_\mathbf{r}g_\mathbf{r}\|_2\\
&\le \|(f-f_n)g_{\mathbf{r}}\|_2+\|(f_n-f)_\mathbf{r}g_\mathbf{r}\|_2\\
&\le \|g_\mathbf{r}\|_{\infty}\|f_n-f\|_2.
\end{align*}
which gives that $(fg)_\mathbf{r}=f_\mathbf{r}g_\mathbf{r}$.
\end{proof}

\begin{theorem}\label{t1.7}
For $f\in L^p(\mathbb{T}^\omega)$, $1\le p<\infty$,  it follows that $f_{\mathbf{r}}\to f$ with respect to $L^p$-norm as $r_1\to 1$. In particular, for $f\in C(\mathbb{T}^\omega)$, we have $f_\mathbf{r}\to f$ uniformly.
\end{theorem}

\begin{proof}
For trigonometric polynomial $h=\sum_{\mathbf{n}\in A}a_\mathbf{n}\gamma_\mathbf{n}$ where $\hat{h}$ supported on a finite set $A\in \mathbb{Z}^\infty$, it follows that $$|h-h_\mathbf{r}|\lesssim \sum_{\mathbf{n}\in A}|a_\mathbf{n}|\left(1-\prod_{k=1}^\infty \left(\frac{1}{p_k}\right)^{\frac{-\log r_1}{\log 2}n_k}\right)\to 0,$$
since the length of elements in $A$ is bounded. By this we have $p_\mathbf{r}\to p$ uniformly. From Proposition \ref{p1.6} we deduce that
\begin{align*}
\|f-f_\mathbf{r}\|_p&\le \|h-h_\mathbf{r}\|_p+\|f-h\|_p+\|f_\mathbf{r}-h_\mathbf{r}\|_p\\
&\le \|h-h_\mathbf{r}\|_p+2\|f-h\|_p,
\end{align*}
where $p\in [1, \infty]$.
Then the density of trigonometric polynomials in $L^p(\mathbb{T}^\omega)$ for $p\in [1, \infty)$ and in $C(\mathbb{T}^\omega)$ completes this proof.
\end{proof}

\begin{corollary}\label{c1.3}
For every $f\in C(\mathbb{T}^\omega)$, if we define $\tilde{f}(\sigma\cdot \boldsymbol{\lambda})=f_\mathbf{r}(\boldsymbol{\lambda}), \mathbf{r}=\mathbf{p}^{-\sigma}$ on $\mathfrak{D}$ and $\tilde{f}(\boldsymbol{\infty})=\hat{f}(\mathbf{0})$, then $\tilde{f}$ is continuous with respect to the compact topology as in Remark \ref{r1.1}.
\end{corollary}

\begin{proof}
Suppose that $f\in C(\mathbb{T}^\omega)$,
Theorem \ref{t1.7} shows that  $\tilde{f}(\sigma\cdot \boldsymbol{\lambda})$ is continuous with respect to the variable $\sigma$ for $\sigma\in [0,\infty)$, and from the definition of $f_\mathbf{r}$, we deduce that $\tilde{f}(\sigma\cdot \boldsymbol{\lambda})$ is continuous with respect to the variable $\boldsymbol{\lambda}$ for $\sigma\in [0,\infty)$. It is sufficient to prove that $\tilde{f}$ is continuous at the point $\boldsymbol{\infty}$. The equality
$$p_{\mathbf{r}}(\boldsymbol{\lambda})=\sum_{\mathbf{n}\in \mathbb{Z}^\infty} \hat{p}(\mathbf {n})\mathbf{p}^{-\sigma\mathbf{n}}\gamma_{\mathbf{n}}(\boldsymbol{\lambda}),$$
which implies that $p_{\mathbf{r}}(\boldsymbol{\lambda})\to \hat{p}(\mathbf{0})$ as $\sigma \to \infty$ holds for every trigonometric polynomial $p$ and so does for each $f\in C(\mathbb{T}^\omega)$.
\end{proof}

\begin{corollary}\label{c1.2}
For $f\in A_\tau(\mathbb{T}^\omega)$ and $\mathbf{r}\in \mathfrak{D}$, it follows that
$$e^f\ast P_{\mathbf{r}}=e^{f\ast P_{\mathbf{r}}}.$$
\end{corollary}
\begin{proof}
From Theorem \ref{t1.5} and Theorem \ref{t1.7}, for a fixed $f\in A_\tau(\mathbb{T}^\omega)$ we have
\begin{align*}
e^f\ast P_{\mathbf{r}}=\left(\sum_{k=0}^\infty\frac{f^k}{k!}\right)\ast P_{\mathbf{r}}=\sum_{k=0}^\infty\frac{f^k\ast P_{\mathbf{r}}}{k!}=e^{f\ast P_{\mathbf{r}}}.
\end{align*}
where the summation is with respect to the uniform norm.
\end{proof}
\section{Main theorem}
This section is devoted to our main theorem, whose bidisk version was applied to prove the Szeg\"{o}'s theorem in the bidisc case \cite{n}. In our proof of Theorem \ref{t1.8}, we circumvented this theorem by employing Theorem \ref{t1.10}, yet it still retains significant theoretical value, which is why we have given it prominent discussion in our paper.

\begin{theorem}
Suppose $w\in L^1(\mathbb{T}^\omega)$ with $\log w\in L^1(\mathbb{T}^\omega)$, then $w=|g|^2$ for some outer function $g$ in $H^2_+(\mathbb{T}^\omega)$ if and only if $\hat{\log w}$ is supported on $\mathbb{Z}^\infty_+\cup (-\mathbb{Z}^\infty_+)$.
\end{theorem}

\begin{proof}
To prove the necessity , suppose that $\supp(\hat{\log w})\in \mathbb{Z}^\infty_+\cup (-\mathbb{Z}^\infty_+)$, Let $\mathbf{r}=\mathbf{p}^{-\sigma}$ with $\sigma> 1$. It follows that $\log w\ast P_\mathbf{r}=\log w\ast  \mathfrak{P}_\mathbf{r}$, where $\mathfrak{P}_\mathbf{r}$ is the measure such that $\hat{\mathfrak{P}_\mathbf{r}}(\mathbf{n})=\mathbf{r}^{|\mathbf{n}|}$ for $\mathbf{n}\in \mathbb{Z}^\infty_+\cup (-\mathbb{Z}^\infty_+)$ and $\hat{\mathfrak{P}_\mathbf{r}}(\mathbf{n})=0$ for other $\mathbf{n}\in \mathbb{Z}^\infty$. The existence of $\mathfrak{P}_\mathbf{r}$ is due to that $\sum_{\mathbf{n}\in \mathbb{Z}^\infty_+}\mathbf{r}^\mathbf{n}=\zeta(\sigma)$, and the Riemann-zeta function $\zeta(\sigma+it)$ converges when $\sigma>1$. Similarly, it can be deduced that there is a continuous function $\mathfrak{C}_\mathbf{r}:=2\sum_{\mathbf{n}\in \mathbb{Z}^\infty_+}\mathbf{r}^\mathbf{n}(\cdot)^\mathbf{n}-1$ on $\mathbb{T}^\omega$ and $\mathfrak{P}_\mathbf{r}=Re \mathfrak{C}_\mathbf{r}$.

For convenience, from now we write $\mathbf{r}=\sigma \cdot \mathbf{1}$ as in Remark \ref{r1.1}.
Note that $\sigma\cdot \mathbf{1}\in \mathfrak{D}$ belongs to $l^1$ when $\sigma \geq 1$,  by this $\log|w|\ast \mathfrak{C}_{2\cdot \mathbf{1}}$ is a continuous function on $\mathbb{T}^\omega$ since its Fourier transform belongs to $l_1(\mathbb{Z}^\infty)$.
There is a function $g'=e^{\frac{1}{2}\log |w|\ast \mathfrak{C}_{2\cdot \mathbf{1}}}$. It is clear that $g'$ is an outer function. Combining this, Theorem \ref{t1.9}\& \ref{t1.10}\& \ref{t1.11} and the theorem in \cite[page 103]{h},  we have $|g'|=|g_{2\cdot \mathbf{1}}|$, where $g$ is an outer function in $H^2_{\tau}(\mathbb{T}^\omega)$ with $w=|g|^2$. then again by Theorem \ref{t1.9}, $g'=\theta g_{2\cdot \mathbf{1}}$ with $\theta$ a constant of norm 1. Since $\supp(\hat{g})=\supp(\hat{g_\mathbf{r}})$, we have that $g\in H^2_+(\mathbb{T}^\omega)$ and $\supp(\hat{\log w})\subset \mathbb{Z}^\infty_+ \cup (-\mathbb{Z}^\infty_+)$.

Conversely, to prove sufficiency, we suppose that:

(1) $w=|g|^2$ for outer function $g\in H^2_+(\mathbb{T}^\omega)$, $\hat{g}(0)=1$.

(2) There exists $\mathbf{n}_0\notin \mathbb{Z}^\infty_+\cup (-\mathbb{Z}^\infty_+)$ with $\mathbf{n}_0\in \supp(\hat{\log w})$.

Without loss of generality, we assume that $n_0>0$ with respect to $\tau$ since $\log w$ is real and use $\sigma_0$ to denote one of the $\sigma$ large enough that $\|1-g_{\sigma\cdot1}\|_\infty<1$. The choice of $\sigma_0$ ensures that $\log g_{\sigma_0\cdot 1}$ is well defined and lies in $A_{\tau}(\mathbb{T}^\omega)$. Since $e^{\gamma_{\mathbf{n}_0}\log g_{\sigma_0\cdot1}}:=\sum_{k=0}^\infty \frac{(\gamma_{\mathbf{n}_0}\log g_{\sigma_0\cdot1})^k}{k!}$ is invertible in $A_\tau(\mathbb{T}^\omega)$, and it is therefore an outer function. Thus, 
$$\log\left|\int_{\mathbb{T}^\omega} e^{\gamma_{\mathbf{n}_0}\log g_{\sigma_0\cdot1}}dm_{\mathbb{T}^\omega}\right|=\int_{\mathbb{T}^\omega}Re(\gamma_{\mathbf{n}_0}\log g_{\sigma_0\cdot1})dm_{\mathbb{T}^\omega}=\frac{1}{2}Re[\hat{\log w_{\sigma_0\cdot1}}(\mathbf{n}_0)].$$

Since $\supp(\hat{g})\subset \mathbb{Z}^\infty_+$, together with proposition \ref{p1.1} we have
\begin{align*}
\log\left|\int_{\mathbb{T}^\omega} e^{\gamma_{\mathbf{n}_0}\log g_{\sigma_0\cdot1}}dm_{\mathbb{T}^\omega}\right|=\log \left|\int_{\mathbb{T}^\omega}\sum_{k=0}^\infty \frac{(\gamma_{\mathbf{n}_0}g_{\sigma_0\cdot1})^k}{k!}dm_{\mathbb{T}^\omega}\right|=\log 1=0,
\end{align*}
which implies $\hat{\log w}(\mathbf{n}_0)$ must be a pure imaginary number.

Let $g_t(x)=g(x-t)$, and $w_t=|g_t|^2$, $g_t$ is also an outer function and $\hat{\log w_t}(\mathbf{n}_0)=\hat{\log w}(\mathbf{n}_0)\gamma_{\mathbf{n}_0}(t)$. Choose a suitable $t$ such that $\hat{\log{w_t}}(\mathbf{n}_0)$ is not pure imaginary, and this leads a contradiction to the existence of $\mathbf{n}_0$.
\end{proof}
    The approach we employed in the preceding proof can be similarly utilized to investigate the outer parts of functions in $H^1_+(\mathbb{T}^\omega)$. Without loss of generality,  suppose $g\in H^1_+(\mathbb{T}^\omega)$ with $\hat{g}(0)=1$.  The key point is that we may take a sufficiently large $\sigma_0$ such that $g_{\sigma_0}\in A(\mathbb{T}^\omega)$ and  $\|1-|g_{\sigma_0}|\|_\infty<1$. Indeed, we have $\||g_{\sigma_0}|-1\|_\infty\le \|\hat{g}\|_\infty(\zeta(\sigma_0)-1)$, and this ensures the existences of  $\sigma_0$. Define $\tilde{g}_{\sigma_0}=e^{2C(\log |g_{\sigma_0}|)-1}$, where $C$ is the analytic contraction operator with respect to $\tau$. $C(\log |g_{\sigma_0}|)$ is well-defined by \cite[Theorem 8.7.2]{rud}. It is easy to see that $|\tilde{g}_{\sigma_0}|=|g_{\sigma_0}|$ and $\tilde{g}_{\sigma_0}$ is an outer function.
    In summary, we arrive at the following conclusion:
    \begin{proposition}
        For $g\in H^1_+(\mathbb{T}^\omega)$ with $\hat{g}(0)=0$, there exist a sufficiently large $\sigma_0$ such that the outer part $\tilde{g_{\sigma_0}}$ of $g_{\sigma_0}$ takes the form:
$$\tilde{g_{\sigma_0}}=e^{2C(\log|g_{\sigma_0}|)-1},$$
where $C$ is the analytic contraction operator on $A(\mathbb{T}^\omega)$.
    \end{proposition}


\section{Poisson measures on right half-plane}
\quad We have mentioned that the Poisson measures are singular after Definition \ref{d5.5}. This makes it difficult to describe the distribution of their mass. In order to more specifically describe some quantitative properties of these measures, we need to consider their action on more familiar mathematical objects. Roughly speaking, we need to 'restrict' these measures to the imaginary axis of the complex domain.  In this section, we use $\mathbb{H}$ to denote the right half-plane , $\overline{\mathbb{H}}$ to denote 
 the set $\mathbb{H}\cup \{\boldsymbol{\infty}\}$, on which the topology is defined similarly as in Remark \ref{r1.1}, just replace $\mathbb{T}^\omega$ in Remark \ref{r1.1} by the imaginary axis $i\mathbb{R}$.

We first clarify the relationship between 
$\mathfrak{D}$ and $\overline{\mathbb{H}}$. Briefly speaking, 
$\overline{\mathbb{H}}$ can be densely embedded into $\mathfrak{D}$ in a proper manner. Under such an embedding, the imaginary axis is densely mapped into 
$\mathbb{T}^\omega$. Specifically, there is a continuous homomorphism $\pi$ from $\mathbb{R}$ to $\mathbb{T}^\omega: t\to \mathbf{p}^{-it}$. By Kronecker's theorem \cite[Theorem 5.1.3]{rud}, the range of $\pi$ is dense in $\mathbb{T}^\omega$. This naturally endows the upper half-plane with the structure of a dense subset of $\mathfrak{D}$, with respect to the topology in Remark \ref{r1.1}, and the corresponding (continuous) embedding map $\pi'$ defines as $$\sigma+it\mapsto \sigma\cdot\pi(t),$$
where $\sigma\in [0, \infty]$, $t\in \mathbb{R}$ and $(\infty,t)$ serves as the point $\boldsymbol{\infty}$.
By the continuity of $\pi$, for $f\in A_\tau(\mathbb{T}^\omega)$ and $\tilde{f}$ as in corollary \ref{c1.3}, the function $\mathfrak{f}(\sigma+it):=\tilde{f}(\pi'(\sigma+it))$ with $\mathfrak{f}(\boldsymbol{\infty}):=\hat{f}(\mathbf{0})$ is continuous on $\overline{\mathbb{H}}$, together with Proposition \ref{p1.6}, we have 
\begin{align}\label{e7.1}
\sup_{\sigma\in [0,\infty], t\in \mathbb{R}}|\mathfrak{f}(\sigma+it)|=\sup_{t\in \mathbb{R}}|\mathfrak{f}(it)|=\sup_{\boldsymbol{\lambda}\in \mathbb{T}^\omega}|f(\boldsymbol{\lambda})|.
\end{align}

From now on, we use $A_\tau(\overline{\mathbb{H}})$ to denote the set of $\mathfrak{f}$ corresponding to $f\in A_\tau(\mathbb{T}^\omega)$ as described above ,  $A_\tau(i\mathbb{R})$ to denote the restriction of the elements of 
$A_\tau(\overline{\mathbb{H}})$ to the imaginary axis $i\mathbb{R}$, $S$ to denote the restriction map, and $\mathfrak{r}$ to denote the map $f\to \mathfrak{f}$. Equation \eqref{e7.1} shows that $A_\tau(i\mathbb{R})$ is isometric to $A_\tau(\mathbb{T}^\omega)$, and so is $A_\tau(\overline{\mathbb{H}})$.

By combining the Riesz representation theorem for positive linear functional on locally compact Hausdorff spaces, the Hahn-Banach extension theorem and equation \eqref{e7.1}, we have that each point $\sigma+it$ in $\overline{\mathbb{H}}$ has a representing measure $\mathfrak{p}_{\sigma+it}$ on $i\mathbb{R}$ such that
\begin{align}\label{e7.2}
\tilde{f}(\mathbf{p}^{-\sigma-it})=\mathfrak{f}(\sigma+it)=\int_{i\mathbb{R}} S(\mathfrak{f}) d\mathfrak{p}_{\sigma+it},\ \forall \mathfrak{f}\in A_\tau(\overline{\mathbb{H}}),
\end{align}
It is obvious that the representing measures of points in $i\mathbb{R}$ are Dirac measures. Thus it just remains to consider the points $\sigma+it$ with $\sigma\in (0, \infty)$. 
Applying \eqref{e7.2} to the characters on $\mathbb{T}^\omega$, we have
\begin{align}\label{e7.3}
\tilde{\gamma}_{\mathbf{n}}(\mathbf{p}^{-\sigma-it})=\mathfrak{r}(\gamma_\mathbf{n})(\sigma+it)=\int_{i\mathbb{R}}S[\mathfrak{r}(\gamma_\mathbf{n})]d\mathfrak{p}_{\sigma+it}.
\end{align}
for all $\mathbf{n}\in \mathbb{Z}^\infty_{\tau^+}$.


Combining \eqref{e7.3} and the positivity of $\mathfrak{p}_{\sigma+it}$,  we obtain
\begin{align}\label{e7.4}
\mathbf{p}^{-|\mathbf{n}|(\sigma+it)}=e^{-(\sigma+it)\tau(|n|)}=\int_{i\mathbb{R}}\mathbf{p}^{-\mathbf{n}x}d\mathfrak{p}_{\sigma+it}(x)=\int_{i\mathbb{R}}\mathbf{p}^{-\mathbf{n}(x-it)}d\mathfrak{p}_{\sigma}(x), \forall \mathbf{n}\in \mathbb{Z}^\infty_{\tau^+} 
\end{align}
From this we have $\mathfrak{p}_{\sigma+it}(E)=\mathfrak{p_\sigma}(E-it)$ for each Borel set on $i\mathbb{R}$, and the Fourier-Steiltjes transform of $\mathfrak{p}_\sigma$ coincides with the function $e^{-\sigma |\cdot|}$ on a dense subset $\{\log q: q\in \mathbf{Q}, q>0\}$ of $\mathbb{R}$. If we further require that $\mathfrak{p}_{\sigma}$ is absolutely continuous, then the choice of such a measure must be unique.

It is clear that $e^{-\sigma|\cdot|}$ is integrable on $i\mathbb{R}$, and its Fourier transform is the function $\frac{\sigma}{\pi(\sigma^2+(\cdot)^2)}$. For every $\sigma>0$, this is also an integrable function on $i\mathbb{R}$. The inversion theorem implies 
that $e^{-\sigma|\cdot|}$ is the Fourier transform of $\frac{\sigma}{\pi(\sigma^2+(\cdot)^2)}$.

Based on the above discussion, we finally obtain the following theorem:
\begin{theorem}\label{t7.1}
For each $\sigma+it$ with $\sigma\in (0,\infty)$, the corresponding representing measure for $A_\tau(\overline{\mathbb{H}})$, that is absolutely continuous on $i\mathbb{R}$, is $\frac{\sigma}{\pi(\sigma^2+(\cdot-it)^2)}dm_{i\mathbb{R}}$. 
\end{theorem}
\begin{remark}\label{r7.2}

(i) For $f\in A_\tau(\mathbb{T}^\omega)$,  the map $f\mapsto S[\mathfrak{r}(f)]$ induces a correspondence between formal Fourier series and formal Dirichlet series. That is 
\begin{align}\label{e7.5}
\sum_{\mathbf{n}\in \mathbb{Z}^\infty_{\tau^+}}a_\mathbf{n}\gamma_\mathbf{n}\mapsto \sum_{q\in \mathbb{Q}, q\geq 1}b_q q^{-it}, b_{e^{\tau(\mathbf{n})}}=a_\mathbf{n},
\end{align}
where $f$ coincides with $\sum_{\mathbf{n}\in \mathbb{Z}^\infty_{\tau^+}}a_\mathbf{n}\gamma_\mathbf{n}$ on $\mathbb{T}^\omega$. Furthermore, for $\tilde{f}$ as in corollary \ref{c1.3}, the map $\tilde{f}\mapsto \mathfrak{f}$ induces that 
\begin{align}\label{e7.6}
  \sum_{\mathbf{n}\in \mathbb{Z}^\infty_{\tau^+}}a_\mathbf{n}\mathbf{r}^\mathbf{n}\gamma_\mathbf{n}\mapsto \sum_{q\in \mathbb{Q}, q\geq 1}b_q q^{-\sigma-it}, b_{e^{\tau(\mathbf{n})}}=a_\mathbf{n}, \mathbf{r}=\mathbf{p}^{-\sigma}.
\end{align}
(ii) By Theorem \ref{t7.1} and
 Corollary \ref{c1.3}, we may rewrite the Haar measure $m_{\mathbb{T}^\omega}$. That is, for each $f\in C(\mathbb{T}^\omega)$,
 \begin{align*}
 \int_{\mathbb{T}^\omega}fdm_{\mathbb{T}^\omega}=\lim_{\sigma\to \infty}\int_{i\mathbb{R}}f(\mathbf{p}^{-i\cdot})\frac{\sigma}{\pi(\sigma^2+(\cdot)^2)}dm_{i\mathbb{R}}.
 \end{align*}
 A more common way to write it is 
  \begin{align*}
\int_{\mathbb{T}^\omega}fdm_{\mathbb{T}^\omega}=\lim_{N\to \infty}\frac{1}{2N}\int_{-N}^Nf(\mathbf{p}^{-i\cdot})dm_{i\mathbb{R}}.
 \end{align*}
\end{remark}
Regarding algebra $A_\tau(\overline{\mathbb{H}})$, not much is known about it, but a brief analysis suggests that its elements are holomorphic in the right half-plane. Indeed, by Proposition \ref{p5.6} and equation \eqref{e7.1}, the imagine of analytic trigonometric polynomials under $\mathfrak{r}$, namely, the Dirichlet polynomials, are dense in $A_\tau(\overline{\mathbb{H}})$ with respect to uniform norm. 

Although studying the algebra $A_\tau(\overline{\mathbb{H}})$ on the right half-plane is a challenging task, the above discussion still help us determine which zero-measure sets on the infinite-dimensional torus can carry the full mass of the Poisson measures $\{P_{\mathbf{r}}\}$ with $\mathbf{r}\in \mathfrak{D}\setminus \{\mathbf{0}\}$. Our finding on this issue is as follow:
\begin{theorem}\label{t7.3}
For each $\mathbf{r}\cdot \boldsymbol{\lambda}$ lies in $\mathfrak{D}\setminus \{\mathbf{0}\}$,
the corresponding Poisson measure $P_{\mathbf{r}\cdot \boldsymbol{\lambda}}$ is concentrated on the coset $\{\boldsymbol{\lambda}\mathbf{p}^{-it}\}_{t\in \mathbb{R}}$. Furthermore, for every Borel set $B$ on $\mathbb{T}^\omega$, we have
$$P_{\sigma\cdot1}(B):=\mathfrak{p}_{\sigma}[{\pi^{-1}(B}\cap \{\mathbf{p}^{it}\}_{t\in \mathbb{{R}}})].$$
\end{theorem}

\begin{proof}
we first note that $\{\mathbf{p}^{-it}\}_{t\in \mathbb{R}}$ is Borel, since the image of any closed interval $[a, b]\subset \mathbb{R}$ under the map $\pi$ is a compact set, and so closed in the Hausdorff group $\mathbb{T}^\omega$.  The aforementioned facts similarly indicate that $\pi$ is Borel-measurable, because every open interval $(a,b)$ is mapped by $\pi$ to a closed set minus two single points.

Conversely, for each Borel set $B$ on $\mathbb{T}^\omega$, by the continuity of $\pi$, $\pi^{-1}(B\cap \{\mathbf{p}^{-it}\}_{t\in \mathbb{R}})$ is also a Borel subset of $\mathbb{R}$. Then, for each point $\sigma+it\in \mathbb{H}$, we may define a set function $\mathfrak{S}_{\sigma+it}$ on the Borel algebra of $\mathbb{T}^\omega$, namely, $\mathfrak{S}_{\sigma+it}(B):=\mathfrak{p}_{\sigma+it}[{\pi^{-1}(B}\cap \{\mathbf{p}^{it}\}_{t\in \mathbb{{R}}})]$.  $\mathfrak{S}_{\sigma+it}$ is a positive measure since $\mathfrak{p}_{\sigma+it}$ is a positive measure, this and equation \eqref{e7.3} gives that $\mathfrak{S}_{\sigma+it}$ is equivalent to the measure $P_{\sigma\cdot \mathbf{p}^{-it}}$:
\begin{align*}
\int_{\mathbb{T}^\omega}\gamma_{\mathbf{n}}d\mathfrak{S}_{\sigma+it}=
\int_{i\mathbb{R}}S[\mathfrak{r}(\gamma_{\mathbf{n}})]d\mathfrak{p}_{\sigma+it}=
\tilde{\gamma}_\mathbf{n}(\mathbf{p}^{-\sigma-it})=\int_{\mathbb{T}^\omega}\gamma_\mathbf{n}dP_{\sigma\cdot \mathbf{p}^{-it}},\ \forall \mathbf{n}\in \mathbb{Z
}^\infty. \end{align*}

Regarding other elements in $\mathfrak{D}\setminus\{0\}$ that not associated with any points in $\mathbb{H}$. Denote such an element by $\mathbf{r}_0\cdot \boldsymbol{\lambda}_0$ as in Remark \ref{r1.1}, then for some $t_0\in \mathbb{R}$, $P_{\mathbf{r}_0\cdot \boldsymbol{\lambda}_0}\ast P_{1\cdot (\mathbf{p}^{-it_0}\boldsymbol{\lambda}_0^{-1})}=P_{\mathbf{r}_0\cdot \mathbf{p}^{-it_0}}$. Through the preceding discussion, the measure in right side of this equation is concentrated on $\{\mathbf{p}^{it}\}_{t\in \mathbb{R}}$. By combining this and $P_{1\cdot (\mathbf{p}^{-it_0}\boldsymbol{\lambda}_0^{-1})}$ is the Dirac measure on $\{\mathbf{p}^{-it_0}\boldsymbol{\lambda}_0^{-1}\}$, we have $P_{\mathbf{r}_0\cdot \boldsymbol{\lambda}_0}$ is concentrated on the coset $\{\boldsymbol{\lambda}_0\mathbf{p}^{-it}\}_{t
\in \mathbb{R}}$.
\end{proof}
Due to the singularity of Poisson measures, the classical Corona Theorem—which states that the maximal ideal space of the disk algebra is dense in the maximal ideal space of $H^\infty(\mathbb{T})$ with respect to weak-star topology, is far more difficult to consider in the context of our study. Moreover, for a function $f$ in $H_\tau^p(\mathbb{T}^\omega)$, we cannot expect $f_{\mathbf{r}}$ to be continuous for $\mathbf{r}\in \mathfrak{D}\setminus \mathbb{T}^\omega$. However, it can be proven that such $f_{\mathbf{r}}$ is continuous with respect to the topology of the real line on almost all of the aforementioned cosets.
\begin{theorem}
For every $\mathbf{r}\in \mathfrak{D}$ and $f\in H_\tau^1(\mathbb{T}^\omega)$, the function $f_\mathbf{r}$  is continuous on the coset $\{\boldsymbol{\lambda}\mathbf{p}^{-it}\}_{t\in\mathbb{R}}$ endowed with the topology of the real line for almost every $\boldsymbol{\lambda}\in \mathbb{T}^\omega$.
\end{theorem}
\begin{proof}
The case of $\mathbf{r}=\mathbf{0}$ is obvious. As for the case of $\mathbf{r}\in \mathfrak{D}\setminus \{0\}$, the proof is similar to that for the average continuity of integrable functions. By Theorem \ref{t7.3}, we obtain that
\begin{align*}
|f_\mathbf{r}(\boldsymbol{\lambda}\mathbf{p}^{ih})-f_\mathbf{r}(\boldsymbol{\lambda})|&\le\int_\mathbb{R}|f(\boldsymbol{\lambda}\mathbf{p}^{ih-it})-f(\boldsymbol{\lambda}\mathbf{p}^{-it})|\frac{\sigma}{\sigma^2+{t}^2}dm_{\mathbb{R}}(t)\\
&\le \int_\mathbb{R}|f(\boldsymbol{\lambda}\mathbf{p}^{ih-it})-g(h-t)|\frac{\sigma}{\sigma^2+{t}^2}dm_{\mathbb{R}}(t)\\
&+\int_\mathbb{R}|g(h-t)-g(t)|\frac{\sigma}{\sigma^2+{t}^2}dm_{\mathbb{R}}(t)\\
&+\int_\mathbb{R}|f(\boldsymbol{\lambda}\mathbf{p}^{-it})-g(t)|\frac{\sigma}{\sigma^2+{t}^2}dm_{\mathbb{R}}(t)\\
&=\uppercase\expandafter{\romannumeral1}+\uppercase\expandafter{\romannumeral2}+\uppercase\expandafter{\romannumeral3},
\end{align*}
where $\mathbf{r}=\sigma\cdot 1$, and $g$ is a function in $C_c(\mathbb{R})$ to be determined.  Proposition \ref{p1.6} gives that
\begin{align*}
\int_{\mathbb{T}^\omega}|f_{\mathbf{r}}|dm_{\mathbb{T}^\omega}\le \|f\|_1,
\end{align*}
which implies $f_\mathbf{r}$ is bounded almost every where on $\mathbb{T}^\omega$. In other word, the equality
\begin{align*}
\int_{\mathbb{T}^\omega}f(\boldsymbol{\lambda}y^{-1})dP_{\mathbf{r}}(y)=\int_{\mathbb{R}}f(\boldsymbol{\lambda}\mathbf{p}^{it})\frac{\sigma}{\sigma^2+t^2}dm_{\mathbb{{R}}}(t)<\infty
\end{align*}
holds for almost every $\boldsymbol{\lambda}$ in $\mathbb{T}^\omega$, which shows that $h(t):=f(\boldsymbol{\lambda}\mathbf{p}^{it})$ lies in $L^1(\mathbb{R}, \frac{\sigma}{\sigma^2+t^2})$ for almost every $\boldsymbol{\lambda}$. Thus, for every $\epsilon>0$ and $\boldsymbol{\lambda}$  with $f_{\mathbf{r}}(\boldsymbol{\lambda})<\infty$, we can choose $g\in C_c(\mathbb{R})$ such that $\|g-h\|_{L^1(\mathbb{R}, \frac{\sigma}{\sigma^2+t^2})}\le \epsilon/4$. This 
immediately gives that $\uppercase\expandafter{\romannumeral3}\le \epsilon/3$. Similarly, we write
\begin{align*}
\uppercase\expandafter{\romannumeral1}=\int_\mathbb{R}|f(\boldsymbol{\lambda}\mathbf{p}^{ih-it})-g(h-t)|\frac{\sigma}{\sigma^2+{(t-h)}^2}\frac{\sigma^2+{(t-h)}^2}{\sigma^2+{t}^2}dm_{\mathbb{R}}(t).
\end{align*}
Since $\frac{\sigma^2+{(t-h)}^2}{\sigma^2+{t}^2}$ uniformly tends to $1$ on  $\mathbb{R}$ as $h\to 0$, we deduce that $\uppercase\expandafter{\romannumeral1}\le \epsilon/3$ for sufficient small $h$. Likewise, for sufficient small $h$, $\uppercase\expandafter{\romannumeral2}\le \epsilon/3$ since $g$ lies in $C_c(\mathbb{R})$. This completes the proof.
\end{proof}

\end{document}